\documentclass[12pt,a4paper]{elsarticle}
\usepackage[utf8]{inputenc}
\usepackage{amsmath}
\usepackage{amsfonts}
\usepackage{amssymb}

\newcommand{\Hw}{\mathcal{H}^\omega}
\newcommand{\R}{\mathcal{R}}
\newcommand{\M}{{\bf E}}
\DeclareMathOperator*{\esssup}{{\rm ess\,sup}}
\DeclareMathOperator*{\essinf}{{\rm ess\,inf}}

\allowdisplaybreaks
\newtheorem{theorem}{Theorem}
\newtheorem{lemma}{Lemma}
\newtheorem{corollary}{Corollary}

\journal{Journal of Mathematical Analysis and Applications}

\begin{document}

\begin{frontmatter}


\title{On optimal recovery of integrals of random processes}


\author{Oleg Kovalenko}
\ead{olegkovalenko90@gmail.com}

\address{Oleg Kovalenko, Department of Mathematics and Mechanics, Oles Honchar Dnipro National University, Gagarina~ave., 72, Dnipro, 49010, Ukraine}

\begin{abstract}
 In this paper we prove a sharp Ostrowski type inequality for random processes of certain classes. This inequality is later applied to a solution of the optimal recovery of the integral $\int_0^1\xi_tdt$, using the random variables $\xi_{\tau_1},\dots, \xi_{\tau_n}$ as an information set, where $\tau_1,\dots, \tau_n$ are random variables. We also consider the problem of the information set optimization.

\end{abstract}


\end{frontmatter}


\section{Introduction}
The problems of optimal recovery of operators and functionals is an important domain of Analysis. The first results in this area were obtained by Kolmogorov in $1940$s, where the best linear methods of recovery for integrals (also known as optimal quadrature formulas) were considered.

The problem of optimal quadrature formulas is a partial case of the following general problem of optimal recovery. Let a metric space $(X,h_X)$ and sets  $Z$, $Y$, and $W\subset Z$ be given. Assume also that mappings $\Lambda \colon Z\to X$ and $I\colon W\to Y$  are given. An arbitrary mapping $\Phi\colon Y\to X$ is called a method of recovery of $\Lambda$ on the class $W$ given the informational mapping $I$. The error of recovery of $\Lambda$ by the method $\Phi$ is 
$$
{\cal E}(\Lambda,W,I,\Phi)={\sup_{z\in W}h_X(\Lambda(z), \Phi(I(z)))}.
$$
The value
\begin{equation}\label{errorOfRecovery}
    {\cal E}(\Lambda,W, I)=\inf_{\Phi} {\cal E}(\Lambda,W,I,\Phi)
\end{equation}
is called the optimal error of recovery of $\Lambda$ on the class $W$ based on the information given by the mapping $I$. The problem of optimal recovery is to find the optimal error of recovery ${\cal E}(\Lambda,W, I)$ and the optimal method of recovery $\Phi$ that delivers the infimum in~\eqref{errorOfRecovery}, if it exists. Another problem of interest is to find the best (among some set of addmissible) information operator $I$.

The problem of optimal recovery is heavily studied. Many results in this area can be found in the monographs~\cite{Osipenko}, \cite{Zhensykbaev}, \cite{Sobolev}, \cite{Nikolskyi} and references therein.

We are interested in recovery of the integral of random processes on a class of processes that are determined by a majorant of their modulus of continuity. We cite several results concerning the problems of optimal recovery on the classes of functions that are defined by restrictions on the modulus of continuity of the functions or their derivatives.

In~1968 Korneichuk~\cite{Korneichuk} considered the problem of optimal quadrature formulae on a multivariate class $H^\omega$ of functions with given majorant of modulus of continuity.

In~1974 Motornyi~\cite{Motornyi} found the optimal quadrature formula $\sum\limits_{k=1}^np_kf(x_k)$ for the integral $\int\limits_{0}^{2\pi}f(x)dx$ on the class $W^{r}H^\omega$ of univariate periodic functions with given a convex majorant $\omega$ for the modulus of continuity of the derivatives $x^{(r)}$, $x\in W^{r}H^\omega$, $r>3$ is odd.

In~1975 Drozhzhina~\cite{Drozzhina} obtained the optimal quadrature formulae for the integral $\int\limits_{0}^1\xi_tdt$ using the information about the random variables $\xi_{t_1},\ldots, \xi_{t_n}$ on the class of random processes such that $\left(\M|\xi_t-\xi_s|^2\right)^{\frac 12}\leq \omega(|t-s|)$.

In~1979 Sukharev~\cite{Sukharev} obtained results on optimal quadrature formulae of the class of multivariate functions $f\colon K\to\mathbb{R}$, $K\subset \mathbb{R}^n$ is a measurable set, such that $|f(u)-f(v)|\leq \rho(u,v)$ for all $u,v\in K$ and $\rho\colon K^2\to\mathbb{R}$ satisfies the following conditions.
\begin{enumerate}
\item For each $v\in K$ the function $\rho(\cdot, v)$ is integrable.
\item $\rho(u,v) = \rho(v,u)$, $u,v\in K$.
\item $\rho(u,v)\geq 0$ and $\rho(u,u) = 0$, $u,v\in K$.
\item $\rho(u,v)+\rho(v,w)\geq \rho(u,w)$, $u,v,w\in K$.
\end{enumerate}

Some other related results can be found in \cite{Babenko76}, \cite{Babenko77}, \cite{Babenko95}, \cite{Chernaya1} and~\cite{Chernaya2}.
In~\cite{Bab} the problem of optimal quadratures for set-valued functions was considered. \cite{BBPS} contains results concerning a very broad generalizations of the classes $H^\omega$. 

In 1938 Ostrowski~\cite{Ostrowski38} proved the following theorem. 

{\bf Theorem.} {\it Let $f\colon [-1,1]\to \mathbb{R}$ be a differentiable function and let for all $t\in (-1,1)$ $|f'(t)|\leq 1$ . Then for all $x\in [-1,1]$ the following inequality holds
$$\left|\frac 12\displaystyle \int\limits_{-1}^1 f(t)dt - f(x)\right|
\leq \frac{1+ x^2}{2}.
$$
The inequality is sharp in the sense that for each fixed $x\in [-1,1]$, the upper bound $\frac{1+ x^2}{2}$ cannot be reduced.
}

Being interesting themselves, the Ostrowski type inequalities proved to be a useful tool for solving different extremal problems; in particular Ostrowski type inequalities can be applied to optimal recovery of integrals of the low-smoothness function classes. For some results concerning Ostrowski type inequalities for random variables see~\cite{Barnett}, \cite{Kumar} and references therein. Overall, the topic of Ostrowski type inequalities has a very broad bibliography, see for example~\cite{Pecharic},~\cite{Dragomir} and~\cite{Dragomir17}.

The goal of the article is to prove a sharp Ostrowski type inequality for a specific class of random processes and to apply it to the problem of optimal recovery of the integral. 

The paper is organized as follows. In Section~\ref{sec::notations} we introduce the notations that will be used throughout the paper. In Section~\ref{sec::ostrowskiInequality} we state and prove a sharp Ostrowski type inequality. Section~\ref{sec::optimalRecovery} is devoted to the problem of integral optimal recovery. In Section~\ref{sec::informationSetOptimization} we consider the problem of information set optimization.

\section{Notations}\label{sec::notations}
Let $\{\Omega,\mathcal{F},P\}$ be a probability space. For a random variable $\eta$, defined on the probability space $\{\Omega,\mathcal{F},P\}$, set $\|\eta\|_\infty:=\esssup\limits_{w\in \Omega}|\eta(w)|$.

For  $a>0$ denote by $\R(a)$ the space of all random variables $\eta$ on the space $\{\Omega,\mathcal{F},P\}$ such that $\eta(w)\in [0,a]$ for all $w\in\Omega$. 

Everywhere below $\omega$ denotes a concave modulus of continuity, i.~e. a continuous non-decreasing concave function $\omega\colon [0,\infty)\to[0,\infty)$ such that $\omega(0) = 0$.

For $a>0$ denote by $H^\omega(a)$ the class of functions $x\colon[0,a]\to\mathbb{R}$ such that $|x(s)-x(t)|\leq \omega(|t-s|)$ for all $t,s\in [0,a]$.

For a fixed $\tau\in\R(a)$ denote by $\Hw_\tau(a)$ the set of all measurable random processes $\xi_t$, $t\in [0,a]$, defined on the  probability space $\{\Omega,\mathcal{F},P\}$, and such that for all $\eta\in\R(a)$ 
\begin{equation}\label{HwDefinition}
\M|\xi_\tau - \xi_\theta|\leq \omega(\|\tau-\theta\|_\infty).
\end{equation}

Set $\Hw(a):=\bigcap\limits_{\tau\in\R(a)}\Hw_\tau(a)$, so that $\Hw(a)$ is the class of measurable processes such that inequality~\eqref{HwDefinition} holds for all $\tau, \theta\in \R(a)$.

In the case, when $a=1$, we write $H^\omega$,  $\Hw$, $\Hw_\tau$ and $\R$ instead of $H^\omega(1)$, $\Hw(1)$, $\Hw_\tau(1)$ and $\R(1)$ respectively.

\section{Ostrowski type inequality}\label{sec::ostrowskiInequality}
The following theorem gives an Ostrowski type inequality for random processes of the class $\Hw(a)$.

\begin{theorem}\label{th::OstrowskiInequality}
Let $a>0$ and  $\tau\in\R(a)$ be given. Set 
$t^*:=\left\|\tau(\cdot)-\frac a2\right\|_\infty.$ Then 
\begin{equation}\label{OstrowskiInequality}
\sup\limits_{\xi\in\Hw(a)}\M\left| \int\limits_{0}^a\xi_tdt - a\cdot \xi_\tau \right| 
= 
\int\limits_{0}^{\frac a2-t^*}\omega(s)ds + \int\limits_{0}^{\frac a2+t^*}\omega(s)ds.
\end{equation}
\end{theorem}

We prove the theorem in the case $a=1$, the general case can be proved similarly. The proof of the theorem is given in the following paragraphs.

\subsection{Some remarks about the proof of Theorem~\ref{th::OstrowskiInequality}}\label{p::remarks}

It is enough to prove~\eqref{OstrowskiInequality} for the case of simple random variable $\tau$. The general case can be obtained using approximation of $\tau$ by simple random variables.

Let $\Omega_1,\dots, \Omega_n\in\mathcal{F}$ be pairwise disjoint sets with positive measures such that $P\left(\bigcup\limits_{k=1}^n\Omega_k\right) = 1$, and $
\tau(w) = \tau_k$ for $w\in\Omega_k$, $k=1,\dots, n$.

For a fixed $k\in\{1,\dots, n\}$ set
$$
\tau^*(w):= 
\begin{cases}
\tau(w),&w\in\Omega\setminus\Omega_k,
\\
1-\tau(w),& w\in\Omega_k.
\end{cases}
$$
Since together with arbitrary $\xi\in\Hw_\tau$ (or $\xi\in\Hw$), the process $\xi^*_t$, $t\in[0,1]$,
$$
\xi^*_t(w):= 
\begin{cases}
\xi_t(w),&w\in\Omega\setminus\Omega_k,
\\
\xi_{1-t}(w),& w\in\Omega_k
\end{cases}
$$
belongs to $\Hw_{\tau^*}$ (to $\Hw$ respectively), and for almost all $w\in\Omega$
$$
\left| \int\limits_{0}^1\xi_tdt - \xi_\tau \right|
=
\left| \int\limits_{0}^1\xi^*_tdt - \xi^*_{\tau^*} \right|,
$$
one has 
$$
\sup\limits_{\xi\in\Hw_\tau}\M\left| \int\limits_{0}^1\xi_tdt - \xi_\tau \right| 
=
 \sup\limits_{\xi\in\Hw_{\tau^*}}\M\left| \int\limits_{0}^1\xi_tdt - \xi_{\tau^*} \right|
$$
and $$
\sup\limits_{\xi\in\Hw}\M\left| \int\limits_{0}^1\xi_tdt - \xi_\tau \right| 
=
 \sup\limits_{\xi\in\Hw}\M\left| \int\limits_{0}^1\xi_tdt - \xi_{\tau^*} \right|.
$$

Hence without loss of generality we may assume that $0\leq \tau_k\leq \frac 12$, $k=1,\dots, n$. Moreover, we can also assume that 
$
\tau_1 \leq \ldots\leq \tau_n. 
$

Under the assumptions above, we have $t^* = \frac 12 - \tau_1$ and the right hand side of~\eqref{OstrowskiInequality} becomes 
\begin{equation}\label{rightSide}
\int\limits_{0}^{\tau_1}\omega(s)ds + \int\limits_{0}^{1-\tau_1}\omega(s)ds.
\end{equation}

\subsection{Estimate from above}
\begin{lemma}\label{l::upperEstimate}
Under the assumptions of Theorem~\ref{th::OstrowskiInequality} and of Paragraph~\ref{p::remarks}, the inequality
\begin{equation}\label{upperEstimate}
\M\left| \int\limits_{0}^1\xi_tdt - \xi_\tau \right| 
\leq 
\int\limits_{0}^{\frac 12-t^*}\omega(s)ds + \int\limits_{0}^{\frac 12+t^*}\omega(s)ds
\end{equation}
holds for all $\xi\in \Hw_\tau$, in particular for all $\xi\in \Hw$.
\end{lemma}
For all $\xi\in\Hw_\tau$ one has
\begin{multline}\label{estimateFromAbove1}
\M\left| \int\limits_{0}^1\xi_tdt - \xi_\tau \right|
=
\sum\limits_{k=1}^n \int\limits_{\Omega_k}\left| \int\limits_{0}^1(\xi_t - \xi_{\tau_k})dt \right|P(dw)
\\ \leq 
\sum\limits_{k=1}^n \int\limits_{\Omega_k}\int\limits_{0}^1\left| \xi_t - \xi_{\tau_k} \right|dtP(dw)
= 
\sum\limits_{k=1}^n \int\limits_{0}^1\int\limits_{\Omega_k}\left| \xi_t - \xi_{\tau_k} \right|P(dw)dt
\\ = 
\sum\limits_{k=1}^n \left(\int\limits_{0}^{\tau_k-\tau_1}\int\limits_{\Omega_k}\left| \xi_t - \xi_{\tau_k} \right|P(dw)dt
+
\int\limits_{\tau_k-\tau_1}^{\tau_k}\int\limits_{\Omega_k}\left| \xi_t - \xi_{\tau_k} \right|P(dw)dt\right.
\\+
\left. \int\limits_{\tau_k}^{\tau_k + 1-\tau_n}\int\limits_{\Omega_k}\left| \xi_t - \xi_{\tau_k} \right|P(dw)dt
+
\int\limits_{\tau_k + 1-\tau_n}^{1}\int\limits_{\Omega_k}\left| \xi_t - \xi_{\tau_k} \right|P(dw)dt\right)
\\ = 
\sum\limits_{k=1}^n \left(
\int\limits_{0}^{\tau_k-\tau_1}\int\limits_{\Omega_k}\left| \xi_t - \xi_{\tau_k} \right|P(dw)dt
+
\int\limits_{\tau_k + 1-\tau_n}^{1}\int\limits_{\Omega_k}\left| \xi_t - \xi_{\tau_k} \right|P(dw)dt\right)
\\+
\int\limits_{0}^{\tau_1}\sum\limits_{k=1}^n\int\limits_{\Omega_k}\left| \xi_{\tau_k-s} - \xi_{\tau_k} \right|P(dw)ds
+
\int\limits_{0}^{1-\tau_n}\sum\limits_{k=1}^n\int\limits_{\Omega_k}\left| \xi_{\tau_k + s} - \xi_{\tau_k}  \right|P(dw)dt
\end{multline}
Setting $\theta_s(w):=\tau(w)-s$,  we obtain that $\|\theta_s-\tau\|_\infty = s$, $s\in[0,\tau_1]$, and hence
\begin{equation}\label{estimateFromAbove2}
\int\limits_{0}^{\tau_1}\sum\limits_{k=1}^n\int\limits_{\Omega_k}\left| \xi_{\tau_k-s} - \xi_{\tau_k} \right|P(dw)ds
=
\int\limits_0^{\tau_1}\M|\xi_{\theta_s}-\xi_\tau|ds
\leq
\int\limits_0^{\tau_1}\omega(s)ds.
\end{equation}
Analogously 
\begin{equation}\label{estimateFromAbove3}
\int\limits_{0}^{1-\tau_n}\sum\limits_{k=1}^n\int\limits_{\Omega_k}\left| \xi_{\tau_k + s} - \xi_{\tau_k}  \right|P(dw)dt
\leq
\int\limits_0^{1-\tau_n}\omega(s)ds.
\end{equation}
Now set 
$$\theta_s(w):=\begin{cases}
\tau_k+s,&s\in[1-\tau_n,1-\tau_k],\\
\tau_k + (1-\tau_1-\tau_k-s),& s\in (1-\tau_k,1-\tau_1],
\end{cases}$$ 
$w\in\Omega_k$, $k=1,\dots, n$. Since $\tau_1+\tau_k\leq 1$, $k=1,\dots, n$, $|\theta_s-\tau|\leq s$ for almost all $w\in\Omega$ and $s\in [1-\tau_n,1-\tau_1]$. Hence 
\begin{multline}\label{estimateFromAbove4}
\sum\limits_{k=1}^n \left(
\int\limits_{0}^{\tau_k-\tau_1}\int\limits_{\Omega_k}\left| \xi_t - \xi_{\tau_k} \right|P(dw)dt
+
\int\limits_{\tau_k + 1-\tau_n}^{1}\int\limits_{\Omega_k}\left| \xi_t - \xi_{\tau_k} \right|P(dw)dt\right)
\\ = 
\sum\limits_{k=1}^n \left(
\int\limits_{1-\tau_k}^{1-\tau_1}\int\limits_{\Omega_k}\left| \xi_{1-\tau_1-s} - \xi_{\tau_k} \right|P(dw)ds
+
\int\limits_{1-\tau_n}^{1-\tau_k}\int\limits_{\Omega_k}\left| \xi_{\tau_k+s} - \xi_{\tau_k} \right|P(dw)ds\right)
\\ = 
\int\limits_{1-\tau_n}^{1-\tau_1}\sum\limits_{k=1}^n \int\limits_{\Omega_k}|\xi_{\theta_s}-\xi_{\tau_k}|P(dw)ds
=
\int\limits_{1-\tau_n}^{1-\tau_1}\M|\xi_{\theta_s}-\xi_{\tau}|ds
\leq \int\limits_{1-\tau_n}^{1-\tau_1}\omega(s)ds.
\end{multline}

Finally, inequalities~\eqref{estimateFromAbove1}--\eqref{estimateFromAbove4}, together with observation~\eqref{rightSide}, give inequality~\eqref{upperEstimate}.

The lemma is proved.
\subsection{A random process generated by a function}
The following lemma gives a way to generate random processes from the class $\Hw$, given a function $x\in H^\omega$.
\begin{lemma}\label{l::generatedProcess}
Let $x\in H^\omega$ and $F\in\mathcal{F}$, $P(F)>0$, be given. Then the process $\xi_t=\xi_t(x,F)$,  $t\in [0,1]$,
$$\xi_t(w) := 
\begin{cases}
	\frac 1 {P(F)}x(t),& w\in F, \\
	0,& w\in \Omega\setminus F
\end{cases}
$$
belongs to $\Hw$ and
\begin{equation}\label{generatedProcessExpectation}
\M\xi_t = x(t).
\end{equation}
\end{lemma}

In order to prove that $\xi_t\in\Hw$, it is enough to show, that the inequality
 $$ \M\left|\xi_{\theta_1}-\xi_{\theta_2}\right|\leq
\omega\left(\|\theta_{1} - \theta_{2}\|_\infty\right)$$ 
holds for arbitrary two simple random variables $\theta_1$ and $\theta_2$.

Assume that the pairwise disjoint measurable sets $F_{i}\subset F$ with positive measures are such that  $\theta_k(w) = \theta_{i}^k\in[0,1]$, $w\in F_i$, $i=1,\dots, n$, $k=1,2$, and $P\left(\bigcup\limits_{i=1}^{n}F_{i}\right) = P(F)$. Taking into account that $\omega$ is a non-decreasing concave function, one has
 \begin{multline}\label{MEstimate}
 \M\left|\xi_{\theta_1}-\xi_{\theta_2}\right| 
= \sum\limits_{i=1}^{n}\frac {P(F_i)}{P(F)}
 \left|x(\theta_i^1) - x(\theta_i^2) \right|
 \leq 
\sum\limits_{i=1}^{n}\frac {P(F_i)}{P(F)}
  \omega\left(\left|\theta_{i}^1 - \theta_{i}^2\right|\right) 
\\ \leq
	\omega\left(  \sum\limits_{i=1}^{n}\frac {P(F_{i})}{P(F)}
\left|\theta_{i}^1 - \theta_{i}^2\right|\right) 
\leq \omega\left(  \max\limits_{i=1,\ldots,n}\left|\theta_{i}^1 - \theta_{i}^2\right|\right)
   \leq
\omega\left(\|\theta_{1} - \theta_{2}\|_\infty\right).
\end{multline}

Equality~\eqref{generatedProcessExpectation} follows from the definition of the process.

\subsection{Estimate from below}
Let the assumptions of Paragraph~\ref{p::remarks} hold. Consider the process $\xi_t^*\in\Hw$, built according to Lemma~\ref{l::generatedProcess} with $F=\Omega_1$ and $x(\cdot):=\omega(|\cdot-\tau_1|)\in H^\omega$. Then 
\begin{multline*}
\M\left| \int\limits_{0}^1\xi_t^*dt - \xi_\tau^* \right| 
= 
\M \int\limits_{0}^1\xi_t^*dt
= 
 \int\limits_{0}^1\M\xi_t^*dt
= 
\int\limits_{0}^1\omega(|t-\tau_1|)dt
\\ =
\int\limits_{0}^{\tau_1}\omega(t)dt +  \int\limits_{0}^{1-\tau_1}\omega(t)dt, 
\end{multline*}
which together with~\eqref{rightSide} gives the estimate  
$$
\sup\limits_{\xi\in\Hw}\M\left| \int\limits_{0}^1\xi_tdt - \xi_\tau \right| 
\geq 
\int\limits_{0}^{\frac 12-t^*}\omega(s)ds + \int\limits_{0}^{\frac 12+t^*}\omega(s)ds.
$$

\section{Optimal recovery of integrals}\label{sec::optimalRecovery}
\subsection{Statement of the problem}
Let $n\in\mathbb{N}$ and random variables $\tau_1,\dots, \tau_n\in\R$ be given.
For an arbitrary function $\varphi\colon\R^n\to\R$, the operator 
$$S = S_n^\varphi(\xi;\tau_1,\dots, \tau_n) = \varphi(\xi_{\tau_1},\ldots, \xi_{\tau_n})$$
is called a method of recovery  of the integral $\int\limits_0^1\xi_tdt$ of the random process $\xi\in\Hw$. The number 
 \begin{equation}\label{def::recoveryError}
 e(\tau_1,\ldots,\tau_n; S):= \sup\limits_{\xi\in\Hw}\M\left|\int\limits_0^1\xi_tdt - S^\varphi_{n}(\xi;\tau_1,\ldots,\tau_n)\right|
 \end{equation}
is called the error of recovery of the method $S$. 

We consider the following problem. For fixed $n\in\mathbb{N}$ and $\tau_1,\dots, \tau_n\in\R$, find the value of the optimal recovery
\begin{equation}\label{def::optimalRecovery}
E(\tau_1,\ldots,\tau_n):=\inf\limits_S e(\tau_1,\ldots,\tau_n; S),
\end{equation}
and the optimal method of recovery $S$, on which the infimum in~\eqref{def::optimalRecovery} is attained.

In~\cite{Drozzhina}, the problem of integral optimal recovery~\eqref{def::optimalRecovery} was considered in the case when $\tau_1,\ldots, \tau_n$ are constants on an analogue of the class $\Hw$. 

\subsection{Main result}
For $t\geq 0$ set 
\begin{equation}\label{Idef}
I(t) :=\int\limits_{0}^{t}\omega(s)ds.
\end{equation}
The following theorem gives a solution to the integral optimal recovery problem in a special case, when $\tau_1,\ldots, \tau_n$ contain ''one degree of randomness''.

\newcommand{\rightSide}{2\sum\limits_{k=1}^{n-1}I\left(\frac {t_{k+1}-t_{k}}2\right) +
I\left(\frac {1-t_n}2-t^*\right) + I\left(\frac {1-t_n}2+t^*\right).}
\begin{theorem}\label{th::main}
Let $n\in\mathbb{N}$, $\tau\in \R$ and the numbers $0=t_1<\ldots<t_n$ be such that $\tau+t_n\leq 1$ almost everywhere. Set $\tau_k :=\tau+t_k$, $k=1,\dots, n$, and $t^*:=\left\|\tau - \frac {1-t_n}2\right\|_\infty$. Then
\begin{multline}\label{optimalRecoveryValue}
E(\tau_1,\ldots,\tau_n) = \rightSide
\end{multline}
The optimal recovery method is $S = \sum\limits_{k=1}^nc_k^*\xi_{\tau_k}$, where $c_1^* = \tau+\frac {t_2-t_1}2$, $c_k^* = \frac{t_{k+1}-t_{k-1}}2$, $k=2,\dots, n-1$ and $c_n^*= 1-\tau - \frac {t_n+t_{n-1}}2$.
\end{theorem}
The proof of this theorem will follow from the results of subsequent paragraphs.
\subsection{Auxiliary result}
\begin{lemma}\label{l::shiftedProcess}
Let $a>0$, $\tau\in \R(a)$ and $b>0$ be such that $\tau + b\leq a$ almost everywhere. For a process $\xi\in \Hw(a)$ set
\begin{equation*}
\zeta_t(w) := 
\begin{cases}
	\xi_t(w)-\xi_{\tau}(w), & 0\leq t\leq\tau(w),\\
	\xi_{t+b}(w)-\xi_{\tau +b}(w), &\tau(w)<t\leq a-b.
\end{cases}
\end{equation*}
Then $\zeta\in \Hw_{\tau}(a-b)$ and  $\zeta_\tau \equiv 0$.
\end{lemma}

Equality $\zeta_\tau \equiv 0$ follows from the definition of the process $\zeta$. For a random variable $\theta\in\R(a-b)$ set 
$$\tilde{\theta}(w) =
\begin{cases}
\theta(w), & \theta(w) \leq \tau(w) \\
\theta(w) + b,&  \theta(w) >\tau(w).
\end{cases}
$$
and 
$$\tilde{\tau}(w) =
\begin{cases}
\tau(w), & \theta(w) \leq\tau(w) \\
\tau(w) + b,&  \theta(w) >\tau(w).
\end{cases}
$$
Then
$$\M|\zeta_\theta - \zeta_\tau| =\M|\zeta_\theta| =\M|\xi_{\tilde{\theta}} - \xi_{\tilde{\tau}}|\leq \omega(\|\tilde{\theta} - \tilde{\tau}\|_\infty) =\omega(\|\theta - \tau\|_\infty).$$
Hence $\zeta\in \Hw_\tau(a-b)$ and the lemma is proved.

\subsection{Estimate from above}
In this paragraph we prove that 
\begin{multline}\label{ORupperEstimate}
E(\tau_1,\ldots,\tau_n) \leq  \rightSide
\end{multline}

Set $\alpha_0 :=0$, $\alpha_k :=\tau + \frac {t_k+t_{k+1}}2$, $k=1,\dots, n-1$, and $\alpha_{n} = 1$. Then $c_k^* = \alpha_k-\alpha_{k-1}$, $k=1,\dots, n$. Hence
\begin{multline}\label{ORupperEstimate1}
E(\tau_1,\ldots,\tau_n)
	 \leq 
 \sup\limits_{\xi\in\Hw}\M
\left|\int\limits_{0}^1\xi_tdt-\sum\limits_{k=1}^nc_k^*\xi_{\tau_k}\right|
	\\ =
\sup\limits_{\xi\in\Hw}\M
\left|\sum\limits_{k=1}^n\int\limits_{\alpha_{k-1}}^{\alpha_k}\left(\xi_t-\xi_{\tau_k}\right)dt\right|
	\leq
\sup\limits_{\xi\in\Hw}
\M\left| \int\limits_{0}^{\alpha_1}\left(\xi_t-\xi_{\tau_1}\right)dt 
+
\int\limits_{\alpha_{n-1}}^{1}\left(\xi_t-\xi_{\tau_n}\right)dt \right|
\\ +
\sup\limits_{\xi\in\Hw}\sum\limits_{k=2}^{n-1}\M\int\limits_{\alpha_{k-1}}^{\alpha_k}\left|\xi_t-\xi_{\tau_k}\right|dt.
\end{multline}
Let $\xi\in \Hw$ and $k\in \{2,\dots, n-1\}$. Then
\begin{multline}\label{ORupperEstimate2}
\M\int\limits_{\alpha_{k-1}}^{\alpha_k}\left|\xi_t-\xi_{\tau_k}\right|dt
	=
\M\int\limits_{\tau+ \frac {t_{k-1}+t_{k}}2}^{\tau + \frac {t_k+t_{k+1}}2}\left|\xi_t-\xi_{\tau_k}\right|dt
	=
\M\int\limits_{\tau_k+ \frac {t_{k-1}-t_{k}}2}^{\tau_k + \frac {t_{k+1}-t_k}2}\left|\xi_t-\xi_{\tau_k}\right|dt
	\\=
\M\int\limits_{\frac {t_{k-1}-t_{k}}2}^{\frac {t_{k+1}-t_k}2}\left|\xi_{\tau_k+t}-\xi_{\tau_k}\right|dt
	=
\int\limits_{\frac {t_{k-1}-t_{k}}2}^{\frac {t_{k+1}-t_k}2}\M\left|\xi_{\tau_k+t}-\xi_{\tau_k}\right|dt
	 \leq 	\int\limits_{\frac {t_{k-1}-t_{k}}2}^{\frac {t_{k+1}-t_k}2}\omega(|t|)dt
	 \\=
\int\limits_{0}^{\frac {t_{k}-t_{k-1}}2}\omega(t)dt 
+ 
\int\limits_{0}^{\frac {t_{k+1}-t_{k}}2}\omega(t)dt.
\end{multline}

For arbitrary $\xi\in\Hw$, 
\begin{multline}\label{ORupperEstimate3}
\M\left| \int\limits_{0}^{\alpha_1}\left(\xi_t-\xi_{\tau_1}\right)dt 
	+
\int\limits_{\alpha_{n-1}}^{1}\left(\xi_t-\xi_{\tau_n}\right)dt \right|
	\leq
\M\left| \int\limits_{\tau}^{\tau+\frac {t_2}2}\left(\xi_t-\xi_{\tau_1}\right)dt \right|
		\\+
\M\left| \int\limits_{0}^{\tau}\left(\xi_t-\xi_{\tau_1}\right)dt  	+
\int\limits_{\tau+t_n}^{1}\left(\xi_t-\xi_{\tau_n}\right)dt \right|
	+
\M\left|\int\limits_{\tau+\frac{t_n+t_{n-1}}2}^{\tau+ t_n}\left(\xi_t-\xi_{\tau_n}\right)dt \right|
		\\ \leq
\M\left| \int\limits_{0}^{\tau}\left(\xi_t-\xi_{\tau_1}\right)dt  	+
\int\limits_{\tau+t_n}^{1}\left(\xi_t-\xi_{\tau_n}\right)dt \right|
	+
	 \int\limits_{0}^{\frac {t_2-t_1}2}\M\left|\xi_{\tau_1 + t}-\xi_{\tau_1}\right|dt 
		\\+
\int\limits_{0}^{\frac{t_n-t_{n-1}}2}\M\left|\xi_{\tau_n -t}-\xi_{\tau_n} \right|dt
	\leq
\M\left| \int\limits_{0}^{\tau}\left(\xi_t-\xi_{\tau_1}\right)dt  	+
\int\limits_{\tau+t_n}^{1}\left(\xi_t-\xi_{\tau_n}\right)dt \right|
		\\ +
	 \int\limits_{0}^{\frac {t_2-t_1}2}\omega(t)dt 
	+
\int\limits_{0}^{\frac{t_n-t_{n-1}}2}\omega(t)dt.
\end{multline}
Due to Lemmas~\ref{l::shiftedProcess} and~\ref{l::upperEstimate},
\begin{multline*}
\M\left| \int\limits_{0}^{\tau}\left(\xi_t-\xi_{\tau_1}\right)dt  	
	+
\int\limits_{\tau+t_n}^{1}\left(\xi_t-\xi_{\tau_n}\right)dt \right|
	\leq \sup\limits_{\substack{\zeta\in\Hw_\tau(1-t_n), \\\zeta_\tau\equiv 0} }\M\left| \int\limits_{0}^{1-t_n}\zeta_tdt\right| 
	\\  \leq 
\int\limits_{0}^{\frac {1-t_n}2-t^*}\omega(s)ds + \int\limits_{0}^{\frac {1-t_n}2+t^*}\omega(s)ds
\end{multline*}

The latter inequality, together with inequalities~\eqref{ORupperEstimate1}, \eqref{ORupperEstimate2} and~\eqref{ORupperEstimate3} give the estimate from above~\eqref{ORupperEstimate}.
\subsection{Estimate from below}
Below we prove that 
\begin{multline}\label{ORestimateFromBelow}
E(\tau_1,\ldots,\tau_n) \geq \rightSide
\end{multline}
It is enough to prove this inequality for the case of simple random variable $\tau$ such that assumptions from Paragraph~\ref{p::remarks} hold.

For each $\varphi\colon\R^n\to\R$, taking into account that the class $\Hw$ is centrally symmetric, one has
\begin{multline*}
\sup\limits_{\xi\in\Hw}\M\left|\int\limits_0^1\xi_tdt - \varphi(\xi_{\tau_1},\ldots,\xi_{\tau_n})\right|
\geq
\sup\limits_{\substack{\xi\in\Hw,\, \xi_{\tau_k} \equiv 0, \\ k=1,\ldots, n}}\M\left|\int\limits_0^1\xi_tdt - \varphi(0)\right|
\\ =
\sup\limits_{\substack{\xi\in\Hw,\, \xi_{\tau_k} \equiv 0, \\ k=1,\ldots, n}}\max\left(\M\left|\int\limits_0^1\xi_tdt - \varphi(0)\right|,\M\left|\int\limits_0^1(-\xi_t)dt - \varphi(0)\right|\right)
\\ \geq 
\frac 12\sup\limits_{\substack{\xi\in\Hw,\, \xi_{\tau_k} \equiv 0, \\ k=1,\ldots, n}}\left(\M\left|\int\limits_0^1\xi_tdt - \varphi(0)\right|+\M\left|\int\limits_0^1\xi_tdt + \varphi(0)\right|\right)
\\ \geq \sup\limits_{\substack{\xi\in\Hw,\, \xi_{\tau_k} \equiv 0, \\ k=1,\ldots, n}}\M\left|\int\limits_0^1\xi_tdt \right| 
=  \sup\limits_{\substack{\xi\in\Hw,\, \xi_{\tau_k} \equiv 0, \\ k=1,\ldots, n}}\M\int\limits_0^1\xi_tdt,
\end{multline*}
hence 
\begin{equation}\label{EFromBelow}
E(\tau_1,\ldots,\tau_n) \geq \sup\limits_{\substack{\xi\in\Hw,\, \xi_{\tau_k} \equiv 0, \\ k=1,\ldots, n}}\int\limits_0^1\M\xi_tdt.
\end{equation}

Set $s_0 := 0$, $s_k := \tau_1 + \frac{t_k+t_{k+1}}2$, $k=1,\dots, n-1$ and $s_n:=1$.
Using Lemma~\ref{l::generatedProcess}, define a random process $\xi_t^*:=\xi_t(\Omega_1,x)$, where

\begin{equation}\label{extremalProcessDef1}
x(t) := \omega(|t-(\tau_1 + t_k)|), t\in [s_{k-1},s_k), k=1,\dots,n.
\end{equation}
From the equivalent definition 
$$
x(t) =\min\limits_{k=1,\dots,n} \omega(|t-(\tau_1 + t_k)|),t\in[0,1],
$$
if follows that $x(t)\in H^\omega$, hence $\xi_t^*\in\Hw$. Moreover, since $x(\tau_1 + t_k) = 0$, $k=1,\dots, n$, one has $\xi_{\tau_k}^*\equiv 0$,  $k=1,\dots, n$, and hence, due to~\eqref{EFromBelow}, 
$$
E(\tau_1,\ldots,\tau_n)\geq \int\limits_{0}^{1}\M\xi^*_tdt = \int\limits_{0}^1x(t)dt.
$$
Evaluating the right hand side of the latter inequality, using representation~\eqref{extremalProcessDef1} and the fact that $t^* = \frac {1-t_n}2 -\tau_1$, we obtain the right hand side of~\eqref{ORestimateFromBelow}.

\section{Optimization of the informational set}\label{sec::informationSetOptimization}
\subsection{Measurement times optimization}
In this section we consider the problem of optimization of the information set $\{\tau_1,\ldots, \tau_n\}$, in order to minimize the error of recovery. We consider the random process $\xi_t$ as some physical quantity and the random variables $\xi_{\tau_k}$ to be the measurements of this quantity at (possibly random) times $\tau_k$, $k=1,\ldots, n$.

It appears, that if the error of recovery is measured by the error for the ''worst'' function, as it is defined in~\eqref{def::recoveryError}, then the possibility to choose time for measurements randomly does not give benefits compared to the case, when the measurements are done at some fixed, non-random times. More precisely, the following statement holds.
\begin{corollary}
Under the assumptions of Theorem~\ref{th::main}, 
$$\inf\limits_{\tau_1,\ldots, \tau_n} E(\tau_1,\ldots,\tau_n) = 2nI\left(\frac{1}{2n}\right).$$
The optimal measurement times are given by $\tau_k = \frac{2k-1}{2n}$, $k=1,\dots, n$.
\end{corollary}
Recall, that  $I(t)=\int\limits_{0}^{t}\omega(s)ds$. Since $\omega$ is non-decreasing, $I(\cdot)$ is a convex function. Then for arbitrary $\alpha_1,\ldots \alpha_{2n}>0$ one has 
$$\sum\limits_{s=1}^{2n}I(\alpha_s)\geq 2n I\left(\frac {1}{2n}\sum\limits_{s=1}^{2n}\alpha_s\right)$$
and the statement of the corollary follows from~\eqref{optimalRecoveryValue}.

Let now the measurements be done by such device, that the first measurement is triggered by some random event (which occurs at the random time $\tau_1$) and each of the rest $n-1$ measurements are done at time $\tau_k = \tau_1 + t_k$, i.~e. in $t_k$ time units after the first measurement, $k=2,\dots, n$. The following statement optimizes the choice of the numbers $t_2,\dots, t_n$, given the information about $\tau_1$.
\begin{theorem}\label{th::fixedTauOptimization}
Let the assumptions of Theorem~\ref{th::main} hold and $m :=\essinf\limits_{w\in \Omega}\tau(w)$, $M :=\esssup\limits_{w\in \Omega}\tau(w)$. 

If 
\begin{equation}\label{largeM}
(2n-1)m+M\geq 1,
\end{equation}
then 
$$
\inf\limits_{t_2,\ldots, t_n} E(\tau_1,\ldots,\tau_n) 
=
(2n-1)I\left(\frac{1-M}{2n-1}\right) + I(M)
$$
and the infimum is attained for  $t_k = \frac{2(k-1)(1-M)}{2n-1}$, $k=2,\dots, n$.

If 
\begin{equation}\label{smallM}
(2n-1)M+m\leq 1,
\end{equation}
then 
$$
\inf\limits_{t_2,\ldots, t_n} E(\tau_1,\ldots,\tau_n) 
=
(2n-1)I\left(\frac{1-m}{2n-1}\right) + I(m)
$$
and the infimum is attained for  $t_k = \frac{2(k-1)(1-m)}{2n-1}$, $k=2,\dots, n$.

Otherwise,
$$
\inf\limits_{t_2,\ldots, t_n} E(\tau_1,\ldots,\tau_n) 
=
(2n-2)I\left(\frac{1-m-M}{2n-2}\right) + I(m) + I(M)
$$
and the infimum is attained for  $t_k = \frac{(k-1)(1-m-M)}{n-1}$, $k=2,\dots, n$.
\end{theorem}

The proof of this statement will be given in subsequent paragraphs.
\subsection{Auxiliary results}
Recall, that the vector $a\in\mathbb{R}^d$ majorizes the vector $b\in\mathbb{R}^d$ (denoted by $a\succ b$), iff $\sum\limits_{i=1}^{k}a_{[i]} \geq \sum\limits_{i=1}^{k}b_{[i]}$, $k=1,\dots, d-1$, and $\sum\limits_{i=1}^{d}a_i = \sum\limits_{i=1}^{d}b_i$, where $a_{[i]}$ and $b_{[i]}$ denote the $i$-th biggest coordinates of the vectors $a$ and $b$ respectively.

Karamata's inequality~\cite{Karamata} states, that for every convex function $f$ and vectors $x,y\in\mathbb{R}^d$ such that $x\succ y$, one has
$$\sum\limits_{k=1}^df(x_k)\geq \sum\limits_{k=1}^df(y_k).
$$
We need the following lemma.

\begin{lemma}\label{l::majorization}
Let $x,y \in \mathbb{R}^d$ be such that $x\succ y$ and $a\in\mathbb{R}$. Then 
\begin{equation}\label{extendedMajorization}
(x_1,\ldots, x_d,a)\succ (y_1,\ldots, y_d,a).
\end{equation}
\end{lemma}

It is well known, see for example~\cite[Theorem 2.1]{Arnold}, that $x\succ y$ if and only if there exists a double stochastic matrix $A$ such that $y=Ax$. Then the matrix 
$B=\begin{pmatrix} 
A & 0\\
0 & 1
\end{pmatrix}
$
is also double stochastic and 
$\begin{pmatrix}
y\\a
\end{pmatrix} = 
B
\begin{pmatrix}
x\\a
\end{pmatrix}
$, hence~\eqref{extendedMajorization} holds. The lemma is proved.

A vector $s = (s_1,\dots, s_n)\in\mathbb{R}^n$ with $s_k\geq 0$, $k=1,\dots, n-1$, $s_n\geq M$, $\sum\limits_{k=1}^ns_k = 1$ will be called admissible.

For an admissible vector $s=(s_1,\dots, s_n)$, set 
$$s^M :=\left(\frac {s_1}2,\frac {s_1}2, \frac {s_2}2, \frac {s_2}2,\ldots, \frac {s_{n-1}}2, \frac {s_{n-1}}2, s_n - M, M\right)$$
and 
$$s^m :=\left(\frac {s_1}2,\frac {s_1}2, \frac {s_2}2, \frac {s_2}2,\ldots, \frac {s_{n-1}}2, \frac {s_{n-1}}2, s_n - m, m\right).$$

The following lemmas will be used during the proof of Theorem~\ref{th::fixedTauOptimization}.
\begin{lemma}\label{l::caseA}
Let inequality~\eqref{largeM} hold and $s\in \mathbb{R}^n$ be admissible. Set 
$$
L:=\left(\underbrace{\frac {1-M}{2n-1},\ldots, \frac {1-M}{2n-1}}_{2n-1}, M\right).
$$
Then $s^M\succ L$. Moreover, if 
\begin{equation}\label{snIsBig}
s_n\geq M+m,
\end{equation} 
then $s^m\succ L$.
\end{lemma}

Note, that for arbitrary set $(\alpha_1,\dots, \alpha_d)\in\mathbb{R}^d$, $\alpha:=\frac 1d\sum\limits_{k=1}^d\alpha_k$, one has $(\alpha_1,\dots, \alpha_d)\succ (\alpha,\alpha,\ldots,\alpha)\in\mathbb{R}^d$. Hence, due to Lemma~\ref{l::majorization}, for arbitrary admissible vector $s$, $s^M\succ L$. 

If~\eqref{snIsBig} holds, then $(s_n-m,m)\succ (s_n-M, M)$ and hence, due to Lemma~\ref{l::majorization}, $s^m\succ s^M\succ L$. The lemma is proved.

\begin{lemma}\label{l::caseB}
Let inequality~\eqref{smallM} hold and $s\in \mathbb{R}^n$ be admissible. Set 
$$
L:=\left(\underbrace{\frac {1-m}{2n-1},\ldots, \frac {1-m}{2n-1}}_{2n-1}, m\right).
$$
Then $s^m\succ L$. If 
\begin{equation}\label{snIsSmall}
s_n\leq M+m,
\end{equation} 
then $s^M\succ L$.
\end{lemma}

The proof is  similar to the proof of Lemma~\ref{l::caseA}.


\begin{lemma}\label{l::caseC}
Let neither of inequalities~\eqref{largeM} and~\eqref{smallM} hold and $s\in \mathbb{R}^n$ be admissible. Set 
$$
L:=\left(\underbrace{\frac {1-m-M}{2n-2},\ldots, \frac {1-m-M}{2n-2}}_{2n-2}, m, M\right).
$$
If inequality~\eqref{snIsBig} holds, then $s^m\succ L$. If inequality~\eqref{snIsSmall} holds, then $s^M\succ L$.
\end{lemma}

From the conditions of the lemma it follows, that $m < \frac {1-M-m}{2n-2}< M$.

Let inequality~\eqref{snIsBig} hold. Then
\begin{gather*}
\left(\frac {s_1}2,\frac {s_1}2, \frac {s_2}2, \frac {s_2}2,\ldots, \frac {s_{n-1}}2, \frac {s_{n-1}}2, s_n - m\right)
 \succ
\left(\underbrace{\frac {1-s_n}{2n-2},\ldots, \frac {1-s_n}{2n-2}}_{2n-1}, s_n-m\right)
\\ \succ 
\left(\underbrace{\frac {1-m-M}{2n-2},\ldots, \frac {1-m-M}{2n-2}}_{2n-2}, M\right),
\end{gather*}
where the first majorization  follows from Lemma~\ref{l::majorization} and the second one follows from the inequalities $s_n-m\geq M > \frac {1-M-m}{2n-2}$. The inequality $s^m\succ L$ now follows from Lemma~\ref{l::majorization}.

The second statement of the lemma follows from similar arguments, using the inequalities $s_n-M\leq m < \frac {1-M-m}{2n-2}$.

The lemma is proved.

\subsection{Proof of Theorem~\ref{th::fixedTauOptimization}}
The vector with coordinates $s_k:=t_{k+1}-t_k$, $k=1,\dots, n-1$, $s_n:=1-t_n$ is admissible. Obviously, the numbers $t_k$, $k=1,\ldots, n$, are uniquely determined by an admissible vector $s\in\mathbb{R}^n$.

 Note, that if~\eqref{snIsBig} holds, then
$\left\|\tau - \frac {s_n}2\right\|_\infty = \frac {s_n}2-m$ and hence, due to Theorem~\ref{th::main},
$$E(\tau_1,\ldots,\tau_n) =2\sum\limits_{k=1}^{n-1}I\left(\frac {s_{k}}2\right) +
I\left(m\right) + I\left(s_n-m\right).$$
In the case, when~\eqref{snIsSmall},
$\left\|\tau - \frac {s_n}2\right\|_\infty =M- \frac {s_n}2$ and hence, due to Theorem~\ref{th::main},
$$E(\tau_1,\ldots,\tau_n) =2\sum\limits_{k=1}^{n-1}I\left(\frac {s_{k}}2\right) +
I\left(M\right) + I\left(s_n-M\right).$$

The inequalities from below for the value of $E(\tau_1,\ldots, \tau_n)$ now follow from Lemmas~\ref{l::caseA}, \ref{l::caseB} and~\ref{l::caseC} and Karamata's inequality. Thus it is sufficient to show, that the estimates from below are attained. 

Let inequality~\eqref{largeM} hold. Then $\frac {1-M}{2n-1}\leq m$ and hence for the set $s_k^* = \frac{2(1-M)}{2n-1}$, $k=1,\dots, n-1$ and $s_n^* = M+\frac{1-M}{2n-1}$, inequality~\eqref{snIsSmall} holds, thus 
$$
E(\tau_1^*,\ldots,\tau_n^*) =
(2n-1)I\left(\frac{1-M}{2n-1}\right) + I(M),$$
where $\tau_k^*$ are determined by the numbers $s_k^*$, $k=1,\dots, n$.

Let inequality~\eqref{smallM} hold. Then $\frac {1-m}{2n-1}\geq M$ and hence for the set $s_k^* = \frac{2(1-m)}{2n-1}$, $k=1,\dots, n-1$ and $s_n^* = m+\frac{1-m}{2n-1}$, inequality~\eqref{snIsBig} holds, thus 
$$
E(\tau_1^*,\ldots,\tau_n^*) =
(2n-1)I\left(\frac{1-m}{2n-1}\right) + I(m),$$
where $\tau_k^*$ are determined by the numbers $s_k^*$, $k=1,\dots, n$.

Finally, let neither of inequalities~\eqref{largeM} and~\eqref{smallM} hold. Then for the set $s_k^* = \frac{(1-m-M)}{n-1}$, $k=1,\dots, n-1$ and $s_n^* = m+M$
$$
E(\tau_1^*,\ldots,\tau_n^*) =
(2n-2)I\left(\frac{1-m-M}{2(n-1)}\right) + I(m) + I(M),$$
where $\tau_k^*$ are determined by the numbers $s_k^*$, $k=1,\dots, n$.

The theorem is proved.

\bibliographystyle{elsarticle-num}
\bibliography{bibliography}

\begin{thebibliography}{10}
\expandafter\ifx\csname url\endcsname\relax
  \def\url#1{\texttt{#1}}\fi
\expandafter\ifx\csname urlprefix\endcsname\relax\def\urlprefix{URL }\fi
\expandafter\ifx\csname href\endcsname\relax
  \def\href#1#2{#2} \def\path#1{#1}\fi

\bibitem{Osipenko}
K.~Y. Osipenko, Optimal Recovery of Analytic Functions, Nova Publishers, New
  York, 2000.

\bibitem{Zhensykbaev}
A.~A. Zhensykbaev, Problems of Recovery of Operators, Institute of computer
  research, Moscow, Izhevsk, 2003, (in Russian).

\bibitem{Sobolev}
S.~L. Sobolev, V.~Vaskevich, The Theory of Cubature Formulas, Vol. 415 of
  Mathematics and Its Applications, Springer Netherlands, 1997.

\bibitem{Nikolskyi}
S.~M. Nikol'skyi, Quadrature Formulas, Nauka, Moscow, 1979, (in Russian).

\bibitem{Korneichuk}
N.~P. Korneichuk, Best cubature formulas for some classes of functions of many
  variables, Mathematical Notes of the Academy of Sciences of the USSR 3 (1968)
  360--367.

\bibitem{Motornyi}
V.~P. Motornyi, On the best quadrature formula of the form
  $\sum_{k=1}^np_kf(x_k)$ for some classes of differentiable periodic
  functions, Mathematics of the USSR-Izvestiya 8~(3) (1974) 591–620.

\bibitem{Drozzhina}
L.~V. Drozhzhina, On quadrature formulas for random processes, Dopovidi Akad.
  Nauk Ukrain. RSR, Ser. A 9 (1975) 775--777, (in Ukrainian).

\bibitem{Sukharev}
A.~G. Sukharev, Optimal numerical integration formulas for some classes of
  functions of several variables, Sov. Math. Dokl. 20 (1979) 472--–475.

\bibitem{Babenko76}
V.~F. Babenko, Asymptotically sharp bounds for the remainder for the best
  quadrature formulas for several classes of functions, Math. Notes 19~(3)
  (1976) 187--193.

\bibitem{Babenko77}
V.~F. Babenko, On the optimal error bound for cubature formulas on certain
  classes of continuous functions, Anal. Math. 3~(1) (1977) 3--9.

\bibitem{Babenko95}
V.~F. Babenko, On optimization of weight quadrature formulas, Ukrainian Math.
  J. 47~(8) (1995) 1157–1168.

\bibitem{Chernaya1}
E.~V. Chernaya, Asymptotically exact estimation of the error of weighted
  cubature formulas optimal in some classes of continuous functions, Ukrainian
  Math. J. 47~(10) (1995) 1606--1618.

\bibitem{Chernaya2}
E.~V. Chernaya, On the optimization of weighted cubature formulae on certain
  classes of continuous functions, East J. Approx. 1 (1995) 47--60.

\bibitem{Bab}
V.~F. Babenko, V.~V. Babenko, M.~V. Polischuk, On the optimal recovery of
  integrals of set-valued functions, Ukrainian Math. J. 67~(9) (2016)
  1306--1315.

\bibitem{BBPS}
V.~Babenko, Y.~Babenko, N.~Parfinovych, D.~Skorokhodov, Optimal recovery of
  integral operators and its applications, Journal of Complexity 35 (2016)
  102--123.

\bibitem{Ostrowski38}
A.~Ostrowski, \"{U}ber die absolut abweichung einer differentienbaren
  funktionen von ihren integralmittelwert, Comment. Math. Hel 10 (1938)
  226--227.

\bibitem{Barnett}
N.~S. Barnett, S.~S. Dragomir, Ostrowski type inequalities for {L}ebesgue
  integral: a survey of recent results, Australian J. Math. Anal. Appl. 14~(1)
  (2017) 1--–287.

\bibitem{Kumar}
P.~Kumar, An inequality of ostrowski type for cumulative distribution
  functions, Kyungpook Math. J. 39~(2) (1999) 303--311.

\bibitem{Pecharic}
D.~S. Mitrinovic, J.~Pecaric, A.~M. Fink, Inequalities Involving Functions and
  Their Integrals and Derivatives, Kluwer Academic Publishers, Dordrecht, 1994.

\bibitem{Dragomir}
S.~S. Dragomir, T.~M. Rassias, Ostrowski Type Inequalities and Applications in
  Numerical Integration, Springer, Dordrecht, 2002.

\bibitem{Dragomir17}
S.~S. Dragomir, Ostrowski type inequalities for {L}ebesgue integral: a survey
  of recent results, Australian J. Math. Anal. Appl. 14~(1) (2017) 1--–287.

\bibitem{Karamata}
J.~Karamata, Sur une in{\'e}galit{\'e} relative aux fonctions convexes, Publ.
  Math. Univ. Belgrade 1 (1932) 145--–148, (in French).

\bibitem{Arnold}
B.~C. Arnold, Majorization and the Lorenz Order: A Brief Introduction, Vol.~43
  of Lecture Notes in Statistics, Springer-Verlag, 1987.

\end{thebibliography}
\end{document}